 \documentclass[sn-mathphys-num]{sn-jnl}

\usepackage{amsmath,bm,upgreek}
\usepackage{subcaption}
\usepackage{tabularx}
\usepackage{enumitem}
\newlist{steps}{enumerate}{1}

\usepackage{graphicx}%
\usepackage{multirow}%
\usepackage{amsmath,amssymb,amsfonts}%
\usepackage{amsthm}%
\usepackage{mathrsfs}%
\usepackage[title]{appendix}%
\usepackage{xcolor}%
\usepackage{textcomp}%
\usepackage{manyfoot}%
\usepackage{booktabs}%
\usepackage{algorithm}%
\usepackage{algorithmicx}%
\usepackage{algpseudocode}%
\usepackage{listings}%

\theoremstyle{thmstyleone}%
\theoremstyle{thmstyletwo}%

\theoremstyle{thmstylethree}%

\begin{document}

\title[Semi Analytical Solution of a Nonlinear Oblique Boundary Value Problem]{Semi Analytical Solution of a Nonlinear Oblique Boundary Value Problem}


\author*[1]{\fnm{Mriganka Shekhar} \sur{Chaki}}\email{chaki014@umn.edu}

\author[2]{\fnm{Maria  C.} \sur{Jorge}}\email{mcj@aries.iimas.unam.mx}
\equalcont{These authors contributed equally to this work.}

\affil[1]{\orgdiv{The Hormel Institute}, \orgname{University of Minnesota}, \orgaddress{\street{801 16th Ave NE}, \city{Austin}, \postcode{55912}, \state{Minnesota}, \country{USA}}}

\affil[2]{\orgdiv{Unidad Acad\'{e}mica del Instituto de Investigaciones en Matem\'{a}ticas Aplicadas y en Sistemas en el Estado de Yucat\'{a}n} (UA-IIMAS-EY), \orgname{Universidad Nacional Autónoma de M\'{e}xico}, \orgaddress{\street{Tablaje Catastral N°6998, Carretera Mérida-Tetiz Km. 4.5}, \city{Mérida}, \postcode{97357}, \state{Yucat\'{a}n}, \country{M\'{e}xico}}}


\abstract{A new numerical method  is developed to approximate the solution of Laplace's equation in the exterior of the sphere with a strongly non linear boundary value of oblique type.  A functional analysis attempt to solve this type of boundary condition is not straight forward since results about existence and uniqueness of solution are still limited. Hence, a semi analytical method is described here to approach a solution. A perturbation solution  around the monopole convert the nonlinear oblique problem into a series of known Neumann problems in the exterior of the sphere. The monopole was chosen since the problem surged when trying to approximate the exterior gravitational field of the earth, predominantly monopolar, from measurement of its intensity on the surface. The corresponding Green’s function representation for the  exterior Neumann problem gives an exact analytic solution for each perturbation step  as an integral on the surface of the sphere. Nevertheless, the boundary conditions become very complicated and require to be approximated  numerically. The perturbation solutions  given by integrals of the Green's function on the sphere are computed  at each perturbation step  using different  subdivisions (icosahedron, non-uniform mixed-element meshing  and dense triangulation)  of the surface integrals    with the help of adaptive quadrature method. We call icosahedron method to the integration on the sphere with an icosahedron mesh using Gauss 5-point or adaptive quadrature, according to the integration parameter. This method was very effective to deal with the singularity of  the Green's function  successfully avoiding inaccuracies on the numerical approximation and is an important contribution of this work. The numerical perturbation scheme is performed for two given exact solutions. The icosahedron method is found to be very precise.  
  The approximations show the desired properties: they get closer to the exact solutions as the perturbation parameter get smaller, show rapid convergence in the exterior  of the unit sphere and converge to zero as the radius grows. 
}




\maketitle

\section{Introduction}


 This work has a geophysical motivation first studied by Gauss around 1829: Given the magnitude or total intensity of the vector field  ( magnetic or gravitational) on the earth's surface (or on satellite orbits),  find the direction of the exterior field by assuming it to be the gradient of a harmonic scalar potential that vanishes at infinity. 
This problem  has huge applications in geophysical surveys.
Accelerometers can measure the gravitational field in steady state where the
intensity is easier to measure than its direction. A similar situation arrives
for geomagnetic case while measuring with scalar magnetometers (also known
as nuclear precession magnetometers) which only measure the magnitude of
the magnetic field. Although, Superconducting Quantum Interface Devices
(SQUID) can measure full magnetic field, such devices have high cost maintenance. Hence, due to such huge geophysical application in scalar devices,
the problem has puzzled researchers for decades and there have been
only few developments towards its solution.

Questions about existence and uniqueness of the solution have been studied by Backus and other authors, i.e. Jorge \cite{Jorge1}, Jorge and Magnanini \cite{Jorge2}, D\'{i}az \textit{et al.} \cite{Diaz} and Sacerdote and Sans\'{o} \cite{Sacerdote}  giving partial answers to these questions.
Backus proved in  \cite{Backus1} the uniqueness for the gravitational case and the non uniqueness (Backus \cite{Backus2}) for the magnetic case.

Jorge  and Magnanini \cite{Jorge2}
worked on  the solution to Backus' problem for the magnetic field and found  a uniqueness of solution if a boundary condition is given on the magnetic equator. On the other hand, Alberto et al. \cite{Alberto} addressed the non-uniqueness of the solution to the Backus' problem for the magnetic field using   angular momentum algebra and the Clebsch-Gordan coefficients. They  provided another  proof of uniqueness  using  finite expansion of the field by  spherical harmonics and gave  a new linear method to compute the magnetic field from the measurements of the field intensity using  Backus series.
D\'{i}az \textit{et al.} \cite{Diaz}
  studied  a non-linear oblique boundary value problem 
by converting the Backus' exterior problem into an interior problem using Kelvin transformation and proving its well-posedness. 
Later, Kaiser \cite{Kaiser} considered the geomagnetic problem where the direction of the gradient is assumed to be known on the boundary rather than its magnitude. He studied a nonlinear boundary value problem in the exterior domain of a sphere
in two and three dimensions for a given direction field and obtained uniqueness for axisymmetric cases. Much later, Glotov \cite{Glotov} 
considered the original Backus' problem with the data expanded by $\frac{\partial}{\partial n} | \nabla v | $ given on the
boundary of the domain, i.e. normal derivative   of the magnitude of the gradient of the solution. Such problem is also used to estimate the number of sources for the related
inverse source problem in the plane. In the thesis of 
Zheng \cite{Zheng}, the previous works on Backus' problem has been summarized  and the Backus' problem with expanded data has also been discussed in detail.  

 Recently,  Kan et al. \cite{Kan} studied geomagnetic case using Backus problem 
near the dipole in fractional Sobolev spaces. 
Later, the same group Kan et al \cite{Kan2} proved an existence result for the interior Backus problem in the Euclidean ball.
Few  numerical approach to Backus' problem can also be found in the literature.

Macák et al. \cite{Macak, Macak2} used an 
iterative approach applied for solving
the nonlinear satellite-fixed geodetic boundary value problem (NSFGBVP) using the finite
element method (FEM). However, their numerical approaches contain over-simplifications and deviate from the original problem. 

It is to be noted that there exist very few works on the original Backus' problem, especially, for the gravitational field. 
Local existence of the solution for the gravitational case was proved in Jorge \cite{Jorge1}
through a higher order perturbation of the solution $v$ in the space of H\"{o}lder continuous functions. The perturbation gave rise to a series of exterior Neumann problems on the sphere   which were calculated exactly using the corresponding Green's function.
The present work provides a robust numerical approach to the solution of the Backus' problem for the gravitational field.
 A numerical approximation of the perturbation solution to the problem is developed based on the article on local existence of the solution  given in Jorge \cite{Jorge1}. Each term of the perturbation series will be calculated numerically by integrating on the sphere and its exact analytical solution will be given by the Green's function for the exterior Neumann problem.

\section{Mathematical formulation and perturbation scheme}

A unit sphere is considered as an approximation of the earth, hence mathematically the problem can be stated as:
\begin{align}
 \triangle v = 0, \quad \mbox{in} \quad S_e, \qquad |\nabla v| = f,\quad \mbox{in} \quad S, \quad v(x) \to 0 \quad  \mbox{as} \quad  x \to \infty,   \label{eq1}
\end{align}

\noindent where $v$ is the (gravitational or magnetic) field  potential, $\,S\,$  is the surface of the unit sphere $r=1$, $\;S_e\;$ the exterior region $r>1$, and $f$ is a positive real continuous function representing the measure of the intensity of the field on the surface of the sphere. The nonlinearity appearing
in the boundary condition makes the problem very complex because it is an oblique boundary condition.

 Let $v$ be the gravitational potential and assume its intensity is known on the surface of the sphere.
It is well known that the gravitational field is mainly radial, the monopole $1/r$ is the dominant term of the field potential, therefore the gravitational field can be considered as the sum of the dominant field plus a much smaller unknown field $u$. Write $v = 1/r + \varepsilon u$ where $\varepsilon$ is a small positive parameter, then the problem consists of finding the function $u$. It should be harmonic in the exterior of the sphere $S_e$, decay to zero at infinity and should satisfy the perturbed boundary condition on the surface of the sphere given by

\begin{align}
    \nabla v \cdot \nabla v = \frac{1}{r^4} + 2 \varepsilon \frac{\partial u}{\partial n} + \varepsilon^2 \nabla u \cdot \nabla u = f.
    \label{eq2}
\end{align}

\noindent
Without loss of generality  $f$ is replaced by $f^2$ as the boundary condition in \eqref{eq1}  to avoid the use of absolute value in the intensity. When evaluated on the boundary $S$,  $r=1$, it suggests to write the intensity $f^2$ as $(1 + \varepsilon h)$  for some real function $h$.  

\noindent
For the first order perturbation for the potential of the gravitational field, a nonlinear oblique boundary value problem is obtained

\begin{equation}
\begin{split}
 &   \Delta v = \Delta \left(\frac{1}{r} + \varepsilon u\right) = 0, \quad \textrm{in} \quad S_e, \\
 &    \nabla v \cdot \nabla v = 1 + 2 \varepsilon \frac{\partial u}{\partial n} + \varepsilon^2 \nabla u \cdot \nabla u = f^2= 1 + \varepsilon h, \qquad \textrm{on} \quad  S.   
 \end{split}
    \label{eq3}
\end{equation}

\noindent
Since $1/r$ is harmonic in $S_e$, the unknown field $u$  satisfies
\begin{align}
    \Delta u = 0,\quad \textrm {in}\quad S_e, \qquad 2 \frac{\partial u}{\partial n} + \varepsilon \nabla u \cdot \nabla u = h, \quad \textrm{on} \quad S,
    \label{eq4}
\end{align}

 \noindent and since the exterior of the sphere is unbounded, $u$ should go to zero when $r$ goes to infinity.

\noindent
Equation \eqref{eq4} is still an oblique problem hence a higher order perturbation is required,  define
\begin{align}
   &v_n = {1\over r} + \varepsilon u =  {1\over r} + \varepsilon u_1 + \varepsilon^2 u_2 + \varepsilon^3 u_3 + \cdots + \varepsilon^n u_n,
   \label{eq5}
   \end{align}

   \noindent where $u_k$ is harmonic in $S_e$, for $k = 1, 2, ..., n$ and $\varepsilon\ll 1$. When inserted $u = u_1 + \varepsilon u_2 + \varepsilon^2 u_3 + \cdot + \varepsilon^{n-1} u_n$\;  in \eqref{eq4},  the following exterior Neumann problems are obtained:
   \begin{align}
 &   O(\varepsilon): \left[\frac{\partial u_1}{\partial n}\right]_S = {1\over 2} h := B_1, \qquad \Delta u_1 = 0, \quad \textrm{in} \quad S_e,\label{eq6} \\
&        O(\varepsilon^2): \left[\frac{\partial u_2}{\partial n}\right]_S = - \frac{1}{2} \left(\nabla u_1 \cdot \nabla u_1\right)_S := B_2,\qquad \Delta u_2 = 0 \quad \textrm{in} \quad S_e,\label{eq7}\\
&    O(\varepsilon^3): \left[\frac{\partial u_3}{\partial n}\right]_S = - \left(\nabla u_1 \cdot \nabla u_2\right)_S := B_3, \qquad \Delta u_3 = 0\quad \textrm{in} \quad S_e,\label{eq8}\\
 & \hspace{2.2cm}\vdots\nonumber\\
 &\vspace{2cm}\nonumber\\
&       O(\varepsilon^n): \left[\frac{\partial u_n}{\partial n}\right]_S = - \frac{1}{2} \sum_{\substack{1\le k,j\le n\\ k+j=m}} \left(\nabla u_k \cdot \nabla u_j\right)_S := B_n,   \nonumber \\
&  \hspace{1cm}   \Delta u_n = 0 \quad \textrm{in} \quad S_e, \qquad m=2, 3,...,n.\label{eq9}
   \end{align}

\begin{figure}[htp]
    \centering
\includegraphics[scale=0.35]{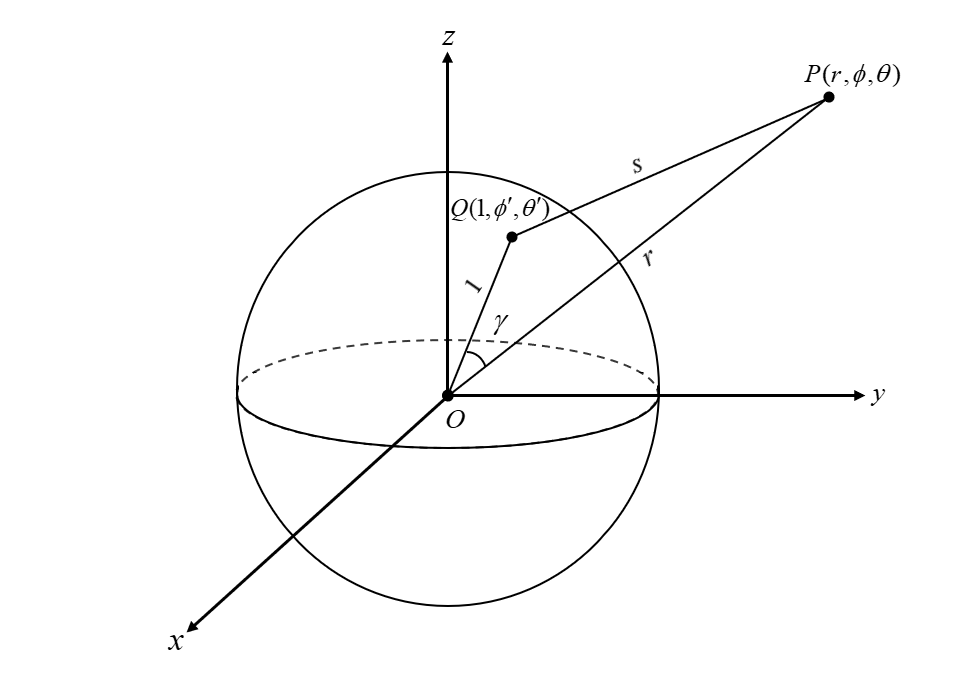}
    \caption{Geometry of Green's function in spherical coordinates, $\phi$ is the longitude and $\theta$ is the co-latitude.}
    \label{fig1}
\end{figure}

\noindent
The Green's function of the second kind for the exterior of the unit sphere is given by \cite{Gubbins}
\begin{align}
    G(P, Q) = \frac{2}{s} - \ln \left [ \frac{1 + s - r \cos \gamma}{r - r \cos \gamma} \right], \label{eq10}
\end{align}

\noindent where \; $\gamma = \measuredangle (OP, OQ)$, \;$\cos \gamma = \cos \theta \cos \theta' + \sin \theta \sin \theta' \cos(\phi - \phi')$,\; $s = \sqrt{1 + r^2 - 2r \cos \gamma}$. The geometry of the Green's function \eqref{eq10} is given in Fig.\ref{fig1} in spherical coordinates where $r$ is the radial coordinate, $\phi$ is longitude and $\theta$ is co-latitude. 


\noindent
Therefore, the solutions for the exterior problems \eqref{eq6}-\eqref{eq9} can be obtained as
 \begin{align}
     u_k(P) = \frac{1}{4\pi} \,\iint_S \, G(P, Q) B_k(Q) dS_Q,\label{eq11}
 \end{align}
 where $ P\in S_e$,\, $ Q\in S$ and $dS_Q$ means integration on the unit sphere. 
 This gives a direct representation for the terms of the perturbation expansion of any order. Hence,  equations \eqref{eq5}-\eqref{eq11}  represent an iterative method to compute an approximation to any order for the gravitational potential. 

Some important properties of Green's function  are given in the  Appendix.

\section{Numerical evaluation}

    In this section, the proposed perturbation scheme is carried out considering a known harmonic function $u$ for the exact solution $v=1/r + \varepsilon u$. 
On following the perturbation scheme,  consider the subsequent steps:

\begin{itemize}

\item[\bf Step 1.] 

Consider first a scalar potential  as the complex harmonic function $u= \frac{1}{r^2}\sin{\theta}\; e^{i\phi}$, to write the exact solution as
\begin{align}
    v(r,\phi,\theta)=\frac{1}{r}+\varepsilon\frac{1}{r^2}\sin{\theta}\; e^{i\phi},
    \label{eq12}
\end{align}
\noindent
where $\varepsilon=10^{-k},  k>0$.

\item[\bf Step 2.]

\begin{figure}[ht]
    \centering
    \includegraphics[scale=0.6]{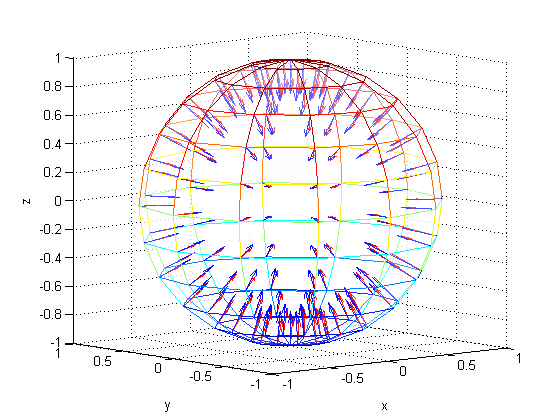}
    \caption{Gravitational field ($\nabla v$) at $12\times 12$ mesh nodes of a unit sphere with red arrows  for  $\varepsilon=10^{-1}$ and blue arrows  for $\varepsilon=10^{-4}$.}
    \label{fig2}
\end{figure}

Next,  based on the considered exact solution  \eqref{eq12},  derive the boundary condition
\begin{align}
   & (\nabla v\cdot \nabla v)_S=\left(\frac{1}{r^4}+\varepsilon\left\{\frac{4}{r^5} \sin{\theta}\; e^{i\phi}+\varepsilon\left(\frac{3}{r^6}\sin^2{\theta} \; e^{2i\phi}\right)\right\}\right)_S,
    \label{eq13}
\end{align}
where
\begin{align}
    \nabla v=\left( -\frac{1}{r^2}-\varepsilon\frac{2}{r^3}\sin{\theta} \; e^{i\phi}\right)\;\hat{r}+\varepsilon\frac{1}{r^3}\cos{\theta}\; e^{i\phi}\;\hat{\theta}+\varepsilon\frac{1}{r^3}  \;ie^{i\phi}\;\hat{\phi}.
    \label{eq14}
\end{align}
Hence, the boundary condition for the unit sphere can be written as
\begin{align}
   (\nabla v\cdot \nabla v)_{r=1}= 1+\varepsilon h(\phi,\theta) = 1 + 4\varepsilon \sin{\theta}\;e^{i\phi}+\varepsilon^2\left( 3\sin^2{\theta}\;e^{2i\phi}\right).\label{eq15}
\end{align}

The gravitational field, i.e. the  gradient vector field   $ \nabla v$  in \eqref{eq14} is represented graphically in Fig. \ref{fig2} by considering a $12\times 12$ surface meshgrid of a unit sphere where red arrows are related to the case of $\varepsilon=10^{-1}$  and blue arrows are of $\varepsilon=10^{-4}$ . Fig. \ref{fig2} shows that as the perturbation parameter gets smaller, the gravitational field vectors become normal to the surface of the unit sphere pointing towards the origin as is expected.

\item[\bf Step 3.]

Calculate $B_1$ using the function $h$ computed in Step 2 for the unit sphere as 
\begin{align}
    B_1=\frac{1}{2}h=2\sin{\theta}\;e^{i\phi}+\varepsilon\left( \frac{3}{2}\sin^2{\theta}\;e^{2i\phi}\right). \label{eq16}
\end{align}

Using $B_1$,    compute  $u_1$   by performing the following surface integration  for the unit sphere 
 \begin{align}
     u_1(P) = \frac{1}{4\pi} \,\iint_S \, G(P, Q) B_1(Q) dS_Q,\label{eq17}
 \end{align}

where the Green's function $G(P,Q)$ is defined in \eqref{eq10} and 
$B_1(Q)$ is defined in \eqref{eq16} for the  exact solution under study.

The above  integration will be computed numerically. Take a  mesh of the surface of the unit sphere and calculate $u_1$ at the nodes of the mesh where  both $P$ and $Q$  are considered on the surface. The value of the integral at some of the mesh points may not be computed because a singularity on the Green's function occurs when $P$ and $Q$ overlaps. Moreover, it can be shown that the Green's function has a singularity when $P$ in $S_e$ and $Q$ on $S$ become collinear with the origen $O$, nevertheless  this singularity is removable.

\item[\bf Step 4.]

In this step,  an approximation function for $u_1$ on the boundary $S$ is constructed by taking a linear combination of $n$-surface harmonics as $u_1(\phi,\theta)=\sum_{k=1}^n \alpha_k\omega_k(\phi,\theta)$, where $\omega_k$'s are  surface harmonics   and $\alpha_k$'s are unknown coefficients. Then consider the meshgrid  on the unit sphere with nodes $P_k, k=1, 2,..., m $ for a given $m$.
Calculating $u_1$ at the nodes $P_k$, a system of $n\times m$ linear equations is obtained which can be solved for the $\alpha_k$'s. Hence, $B_2$ is computed as follows
\begin{align}
\begin{split}
 &   B_2(Q)=-\frac{1}{2}(\nabla u_1\cdot \nabla u_1)_S =- \frac{1}{8} h^2(Q)-\frac{1}{2}\left[ \sum_{k=1}^n \alpha_k \frac{\partial \omega_k (Q)}{\partial \theta}  \right]^2\\
 & \hspace{4.5cm} - \frac{1}{2\sin^2{\theta}} \left[ \sum_{k=1}^n \alpha_k \frac{\partial \omega_k (Q)}{\partial \phi}  \right]^2,
 \end{split}
    \label{eq18}
\end{align}
where
$\nabla u_1(Q)=-\frac{1}{2}h(Q)\;\hat{r} + \sum_{k=1}^n \alpha_k \frac{\partial \omega_k (Q)}{\partial \theta} \;\hat{\theta}+\frac{1}{\sin{\theta}}   \sum_{k=1}^n \alpha_k \frac{\partial \omega_k (Q)}{\partial \phi} \;\hat{\phi}$.
\item[\bf Step 5.]

Compute $u_2$  using  $B_2$ obtained in Step 4 and performing the following surface integration numerically
 \begin{align}
     u_2(P) = \frac{1}{4\pi} \,\iint_S \, G(P, Q) B_2(Q) dS_Q.\label{eq19}
 \end{align}
 The computed $u_2$ will provide the required second approximation $v_2$.

\item[\bf Step 6.]
To compute $v_3$,  construct an approximation function for $u_2$ on the boundary $S$, in a similar fashion as in Step 4, by  taking a  linear combination of $n$-surface harmonics as $u_2(\phi,\theta)=\sum_{k=1}^n \alpha_k^\prime \omega_k^\prime(\phi,\theta)$, where  $\omega_k^\prime$'s are surface harmonics and $\alpha_k^\prime$'s are unknown coefficients. Hence,  compute
\begin{align}
  &  B_3(Q)=- (\nabla u_1\cdot \nabla u_2)_S \nonumber\\
  &\hspace{1cm}=- \frac{1}{2} h(Q) B_2(Q)-\left( \sum_{k=1}^n \alpha_k \frac{\partial \omega_k (Q)}{\partial \theta}  \right)\left( \sum_{k=1}^n \alpha_k^\prime \frac{\partial \omega_k^\prime (Q)}{\partial \theta}  \right) \nonumber\\
  &\hspace{1.5cm} -\frac{1}{\sin^2{\theta}} \left( \sum_{k=1}^n \alpha_k  \frac{\partial \omega_k (Q)}{\partial \phi}  \right)\left( \sum_{k=1}^n \alpha_k^\prime \frac{\partial \omega_k^\prime (Q)}{\partial \phi}  \right),
    \label{eq20}
\end{align}
where 
\begin{displaymath}\nabla u_2(Q)=- B_2(Q)\;\hat{r} + \sum_{k=1}^n \alpha_k^\prime \frac{\partial \omega_k^\prime (Q)}{\partial \theta} \;\hat{\theta}+\frac{1}{\sin{\theta}}   \sum_{k=1}^n \alpha_k^\prime \frac{\partial \omega_k^\prime (Q)}{\partial \phi} \;\hat{\phi}.
\end{displaymath}

\item[\bf Step 7.]

Now, $u_3$ is obtained  using the above $B_3$  and performing the following surface integration 
 \begin{align}
     u_3(P) = \frac{1}{4\pi} \,\iint_S \, G(P, Q) B_3(Q) dS_Q.\label{eq21}
 \end{align}
Finally, the computed $u_3$ will provide the required third approximation $v_3$.

\end{itemize}

 \subsection{Numerical computation of approximation solutions  for the unit sphere }

In order to find the first, second and third approximations  ($v_1$, $v_2$ and $v_3$) of the gravitational potential and corresponding gravitation field approximations ($\nabla v_1$, $\nabla v_2$ and $\nabla v_3$), the numerical integrations on the unit sphere  \eqref{eq17}, \eqref{eq19} and \eqref{eq21} for $u_1$, $u_2$ and $u_3$, respectively,  need to be performed. 
Atkinson \cite{Atkinson} provided a detailed discussion on the numerical integration on the surface of a sphere   where finite element integration using spherical meshing and polyhedron (tetrahedron, octahedron, and  icosahedron) have been discussed. 
On following  \cite{Atkinson},
the numerical integrations \eqref{eq17}, \eqref{eq19} and \eqref{eq21} will be performed by meshing the surface of the unit sphere and computing the values of $u_k$'s and thus $v_k$'s  at the nodes. The following cases are now considered.

\subsubsection{Regular icosahedron}

\begin{figure}[htp]
    \centering
    \begin{minipage}{0.4\textwidth}
\includegraphics[scale=0.4]{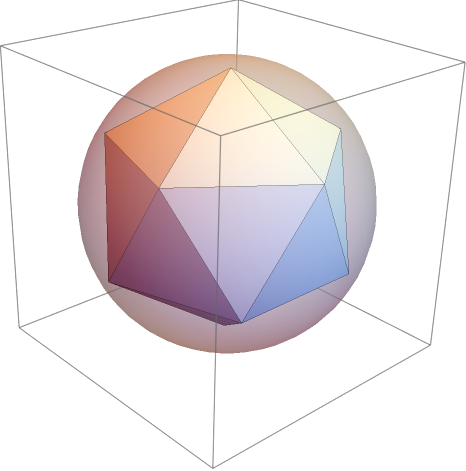}
    \subcaption*{(a)}
    \end{minipage}
     \begin{minipage}{0.4\textwidth}
    \includegraphics[scale=0.5]{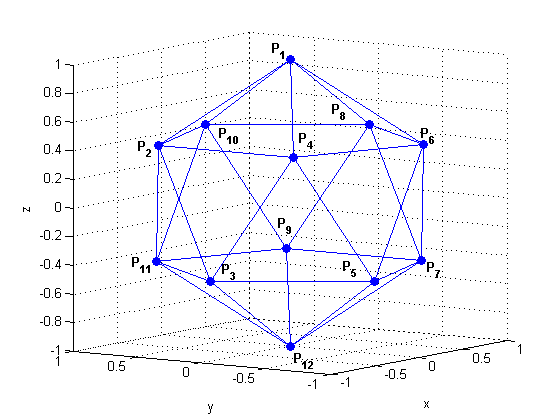}
    \subcaption*{(b)}
     \end{minipage}
    \caption{ (a) An icosahedron inscribed in a unit sphere and (b) MATLAB plot of the icosahedron with 12 vertices. }
    \label{fig3}
\end{figure}

 Construct a uniform mesh of 12 points on the unit sphere by inscribing in $S$ a regular icosahedron which has  20 equal triangles, projecting it on the surface and taking its 12 vertices as the nodes that will constitute the regular mesh on the surface of the unit sphere (see Fig. \ref{fig3}(a)). 
The location of the 12 vertices, i.e.  $P_k(1,\phi_k,\theta_k), k=1,2,...,12$ in  spherical coordinates are given in Table \ref{tab1}, where $\theta_u=\cos^{-1}\left( \frac{\cos{\frac{2\pi}{5}}}{\cos{\frac{2\pi}{5}}-1}\right)$ and
 $\theta_l=\cos^{-1}\left( \frac{\cos{\frac{2\pi}{5}}}{1-\cos{\frac{2\pi}{5}}}\right)$.
A MATLAB plot showing the vertices and edges of the icosahedron is provided in Fig. \ref{fig3}(b). 

\begin{table}
    \centering
    \begin{tabular}{|c|c|c|c|}
    \hline 
  {  $\boldsymbol{P_k}$ }  & $\boldsymbol{P_k(r,\phi_k,\theta_k)}$ & $\boldsymbol{P_k}$   & $\boldsymbol{P_k(r,\phi_k,\theta_k)}$  \\
    \hline 
    \hline
    $P_1$ & $P_1(1,0,0)$ &  $P_7$ & $P_7(1,\frac{6\pi}{5},\theta_l)$\\
      \hline
    $P_2$ & $P_2(1,\frac{\pi}{5},\theta_u)$ & $P_8$ & $P_8(1,\frac{7\pi}{5},\theta_u)$\\
  \hline
  $P_3$ & $P_3(1,\frac{2\pi}{5},\theta_l)$ & $P_9$ & $P_9(1,\frac{8\pi}{5},\theta_l)$\\
    \hline
    $P_4$ & $P_4(1,\frac{3\pi}{5},\theta_u)$ & $P_{10}$ & $P_{10}(1, \frac{9\pi}{5},\theta_u)$\\
      \hline
      $P_5$ & $P_5(1,\frac{4\pi}{5},\theta_l)$ & $P_{11}$ & $P_{11}(1,2\pi,\theta_l)$\\
      \hline
          $P_6$ & $P_6(1,\pi,\theta_u)$ & $P_{12}$ & $P_{12}(1,0,\pi)$\\
      \hline
      \end{tabular}
    \caption{Vertices of icosahedron   in spherical coordinates.}
    \label{tab1}
\end{table}

By performing  numerical integration at the 12 vertices of the icosahedron, the first, second and third approximation of the gravitational potential ($v_1$, $v_2$ and $v_3$)  are  computed on following the perturbation scheme. The computed numerical values 
are provided in the Table 1 and 2 of the Appendix where a comparison between the first approximation $v_1$, second approximation $v_2$ and third approximation $v_3$  with the exact solution $v$ is provided for $\varepsilon=10^{-m},  m=2,4$.
Here, for the computation of $B_2$ in Step 4, an approximation function for $u_1$ is used by considering it as a linear combination of the first 8 spherical harmonics which are provided in  the Appendix.
A more detailed discussion on spherical harmonics can be found in the book of Sternberg and Smith \cite{Sternberg}.
As a result,  an over-determined system with 12 linear equations and 8 unknowns is obtained which can be directly solved in MATLAB.

\begin{figure}
    \centering
       \begin{minipage}{0.4\textwidth}
 \includegraphics[scale=0.38]{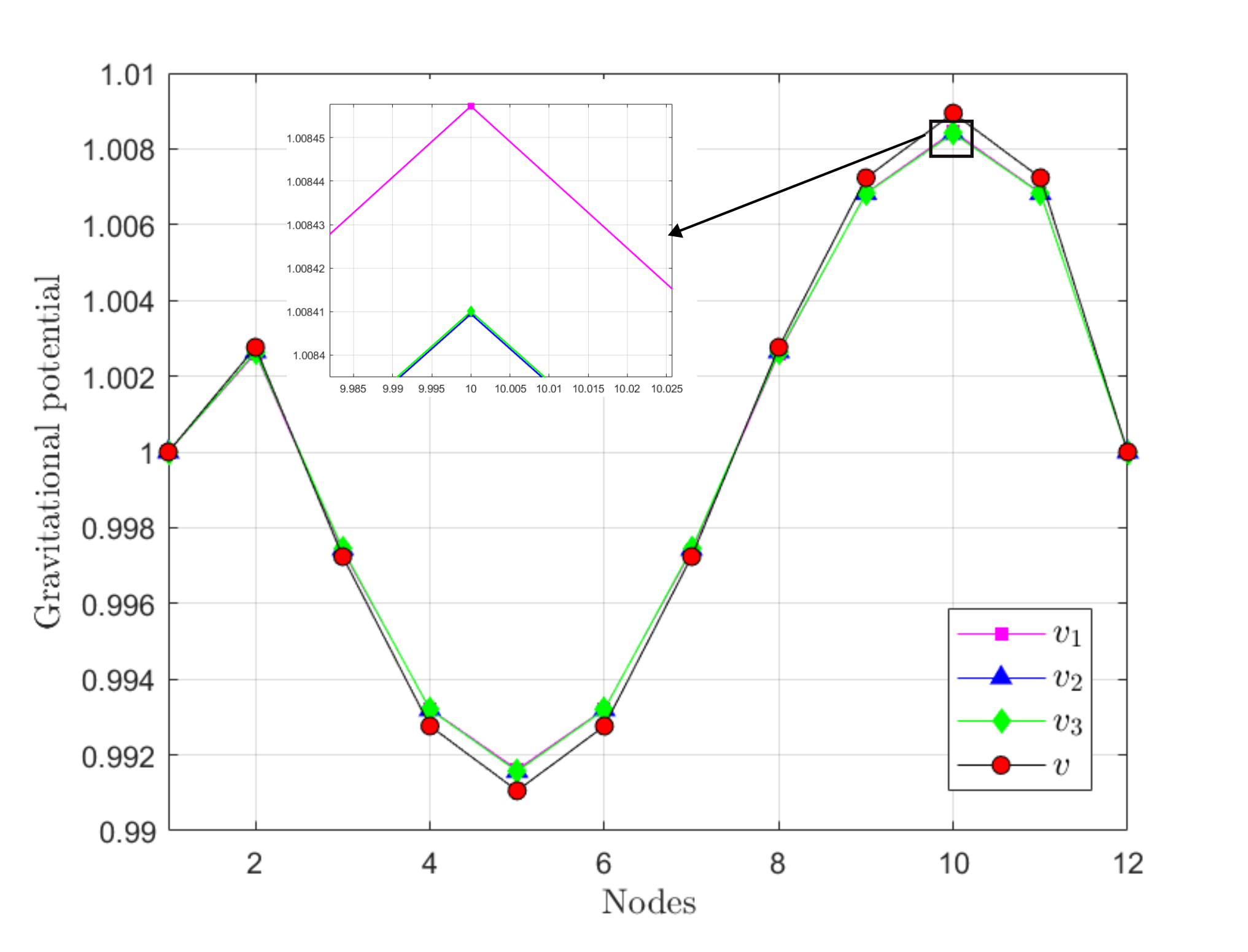}
    \subcaption*{(a)}
    \end{minipage}
    \hspace{1.2cm}
     \begin{minipage}{0.4\textwidth}
    \includegraphics[scale=0.35]{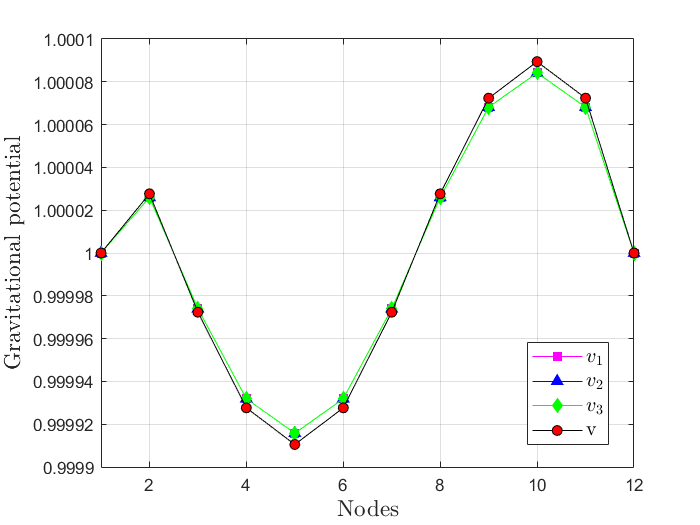}
    \subcaption*{(b)}
     \end{minipage}
    \caption{Graphical representation of the first and second approximation and exact solution against vertices of the icosahedron for (a) $\varepsilon =10^{-2}$ and (b) $\varepsilon =10^{-4}$ for unit sphere.}
    \label{fig4}
\end{figure}

\begin{sidewaystable}
\centering
\fontsize{7.5}{8}
\selectfont
\begin{tabular}{|c|c|c|c|}
\hline
$P$   & $(v-v_1) \times 10^{-5}$                            & $(v-v_2) \times 10^{-5}$                            & $(v-v_3) \times 10^{-5}$                            \\ \hline
$P_1$  & 0.000000000000000 - 0.000000000000000i  & 0.000000000000000 + 0.000000000002071i  & 0.000000000000000 + 0.000000000002072i  \\ \hline
$P_2$  & 0.161760813233158 - 0.496755613709122i  & 0.161377702578847 - 0.497033908187562i  & 0.161377750007574 - 0.497033942661682i  \\ \hline
$P_3$  & -0.161186200620111 - 0.497173094206143i & -0.161569240109127 - 0.496894748024546i & -0.161569287548957 - 0.496894782482404i \\ \hline
$P_4$  & -0.422852870662904 - 0.307478610759380i & -0.422706535418982 - 0.307028320824447i & -0.422706553526719 - 0.307028265062331i \\ \hline
$P_5$  & -0.522894375099003 + 0.000000000000000i & -0.522420912241195 + 0.000000000001438i & -0.522420853610317 + 0.000000000001439i \\ \hline
$P_6$  & -0.422852870662904 + 0.307478610759379i & -0.422706535418982 + 0.307028320829829i & -0.422706553526719 + 0.307028265067715i \\ \hline
$P_7$  & -0.161186200620111 + 0.497173094206137i & -0.161569240109127 + 0.496894748026044i & -0.161569287548957 + 0.496894782483901i \\ \hline
$P_8$  & 0.161760813233158 + 0.496755613709128i  & 0.161377702578847 + 0.497033908191434i  & 0.161377750007574 + 0.497033942665554i  \\ \hline
$P_9$  & 0.422633388175164 + 0.306803113129086i  & 0.422779696251929 + 0.307253486712103i  & 0.422779714370769 + 0.307253430939937i  \\ \hline
$P_{10}$ & 0.522184114859492 + 0.000000000000000i  & 0.522657665658066 + 0.000000000001456i  & 0.522657607016086 + 0.000000000001457i  \\ \hline
$P_{11}$ & 0.422633388175164 - 0.306803113129081i  & 0.422779696251929 - 0.307253486711002i  & 0.422779714370769 - 0.307253430938835i  \\ \hline
$P_{12}$ & 0.000000000000000 - 0.000000000000000i  & 0.000000000000000 + 0.000000000000580i  & 0.000000000000000 + 0.000000000000580i  \\ \hline
\end{tabular}
\caption{Error analysis for  the approximation solutions with exact solution at the vertices of the icosahedron for $\varepsilon=10^{-4}$.}
\label{tab2}
\end{sidewaystable}

In Figs. \ref{fig4}(a) and \ref{fig4}(b), a graphical comparison between the approximation solutions ($v_1$, $v_2$ and $v_3$) and exact solution ($v$) 
is provided for    $\varepsilon=10^{-2}$ and $\varepsilon=10^{-4}$ respectively.  
Fig. \ref{fig4} shows that 
as the perturbation parameter ($\varepsilon$) gets smaller, the approximated solutions ($v_1$, $v_2$ and $v_3$) become close to the exact solution $v$ and 
the higher order approximation  is better  and much closer to the exact solution. An error analysis is provided in Table \ref{tab2} in order to see the convergence of higher order approximations.

The aforementioned numerical integration is not a straight-forward method and it requires some attention. 
 First, it has to be noted that the numerical integration  for the surface integrals were computed  as follows
\begin{align}
\begin{split}
    u_k(P) & =\frac{1}{4\pi}\int_0^\pi \int_0^{2\pi} G(P,\phi^\prime,\theta^\prime) B_k(\phi^\prime,\theta^\prime) d\phi^\prime \sin{\theta^\prime} d\theta^\prime \\
  & =\frac{1}{4\pi}\int_0^\pi I_1(\theta^\prime) \sin{\theta^\prime} d\theta^\prime=\frac{1}{4\pi}\int_{-1}^1 I_1(\zeta)  d\zeta=\frac{1}{4\pi} I_2,
    \end{split}
    \label{eq22}
\end{align}

 \noindent
 for each fixed set of vertices of an inscribed icosahedron $P(1,\phi,\theta)$ where the  change of variable $\zeta=\cos{\theta^\prime}$ was considered. 

 On directly applying the in-built   MATLAB command \verb+`integral2'+ (double integration), i.e. applying the adaptive quadrature method  for double integration with respect to $\phi^\prime$ and $\theta^\prime$ to compute $u_k$, the results can not obtained. This is due to 
the delicate important feature, i.e. the Green's function $G(P,Q)$ has a singularity when $O, P$ and $Q$ becomes collinear and, in particular, on the boundary of the unit sphere, where the points $P$
 and $Q$ exactly overlaps. 
 However, the surface integration is rotational invariant. Hence, by changing the location of the $P$ vertices of the icosahedron 
through rotation, this singularity situation can be avoided and the approximation results can be obtained. 

 Now, another approach to compute $u_k$ is to compute  $I_2$ using Gauss 5-point quadrature formula and  compute $I_1$ using adaptive quadrature along with predetermined error tolerances  which is in-built in the MATLAB command \verb+`integral'+.
 As a result, the integration 
$I_1$ with respect to $\phi^\prime$  is computed for several quadrature points due to sub-interval refinement in adaptive quadrature method whereas  the integration $I_2$ with respect to $\zeta$ (= $\cos{\theta^\prime}$)   is computed for 5  quadrature points. Although, the Green function
has  singularity that occurs whenever $P=Q$,
 this difficulty was overcome on applying two different quadrature methods for $I_1$ and $I_2$ as in \eqref{eq22},  the quadrature points  $Q_k$'s, i.e. ($\phi^\prime, \theta^\prime$) do not overlap with the  fixed vertices of the icosahedron $P_k$'s and hence, it works efficiently providing all the approximations of the gravitational potential given by the perturbation scheme. 

It should be observed that  the results obtained by directly applying in-built   MATLAB command \verb+`integral2'+ for the surface integral along with the rotations to avoid singularity contain more error  as compared to the one obtained using \eqref{eq22}. Hence, the numerical integration method 
explained in \eqref{eq22} by applying two different quadrature methods for $I_1$ and $I_2$ to compute $u_k$ is adopted for further analysis.

\begin{figure}[htp]
    \centering
    \includegraphics[scale=0.6]{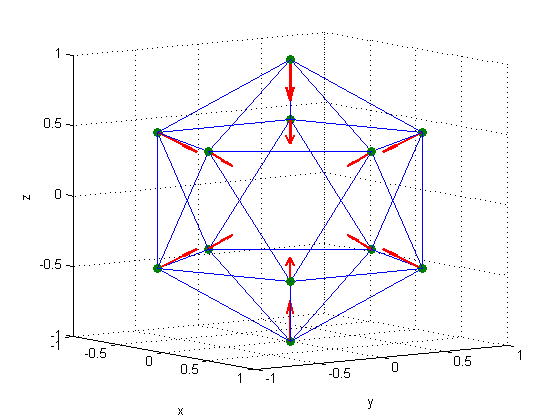}
    \caption{Gravitational field  ($\nabla v_3$)   at 12 vertices of icosahedron for $\varepsilon=10^{-4}$.}
    \label{fig5}
\end{figure}

In Fig. \ref{fig5}, the third approximation of the gravitational field $\nabla v_3$ is drawn showing the vector field at the 12 vertices of the icosahedron for $\varepsilon=10^{-4}$. The similarity between the gravitational field $\nabla v$ of the exact solution shown in Fig. \ref{fig2} and the gravitational field $\nabla v_3$ of the third approximation solution shown in Fig. \ref{fig5} can be observed, this validates the numerical results.

The above results show that an icosahedron is an efficient way of
finding the approximate solution (with a  careful choice of the vertices location  by avoiding  singularity) since
it allows to compute the surface integration much faster. 

\subsubsection{Dense mesh}

\begin{figure}[htp]
    \centering
       \begin{minipage}{0.4\textwidth}
 \includegraphics[scale=0.45]{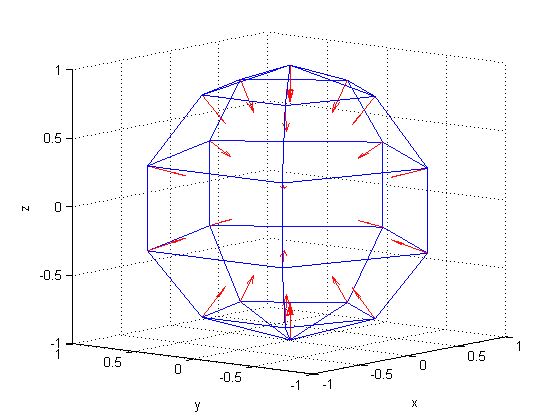}
    \subcaption*{(a)}
    \end{minipage}
    \hspace{1.2cm}
     \begin{minipage}{0.4\textwidth}
    \includegraphics[scale=0.35]{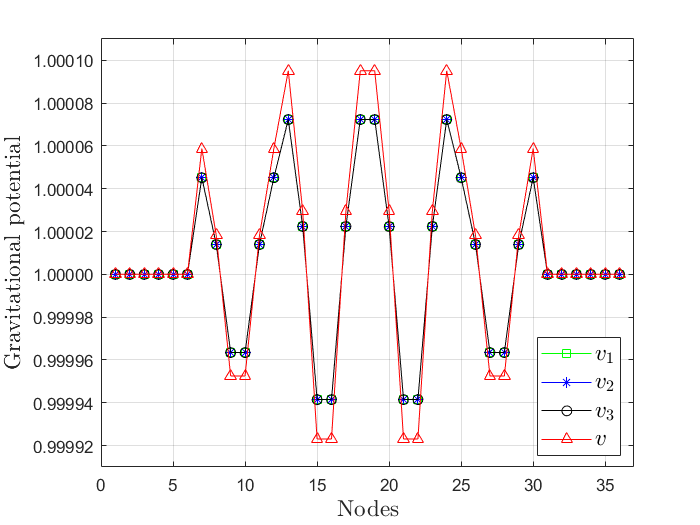}
    \subcaption*{(b)}
     \end{minipage}
    \caption{ (a) Gravitational field approximation $\nabla v_3$ at  36 nodes of mixed (quad/triangular) element meshing of the surface of unit sphere  and (b) comparison of the first, second and third approximation of gravitational potential with exact solution at 36 mesh nodes for  $\varepsilon =10^{-4}$ for the unit sphere.}
    \label{fig6}
\end{figure}

\begin{figure}[htp]
    \centering
       \begin{minipage}{0.4\textwidth}
 \includegraphics[scale=0.45]{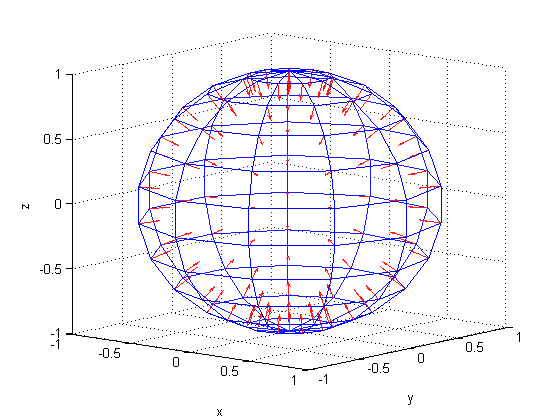}
    \subcaption*{(a)}
    \end{minipage}
    \hspace{1.2cm}
     \begin{minipage}{0.4\textwidth}
    \includegraphics[scale=0.35]{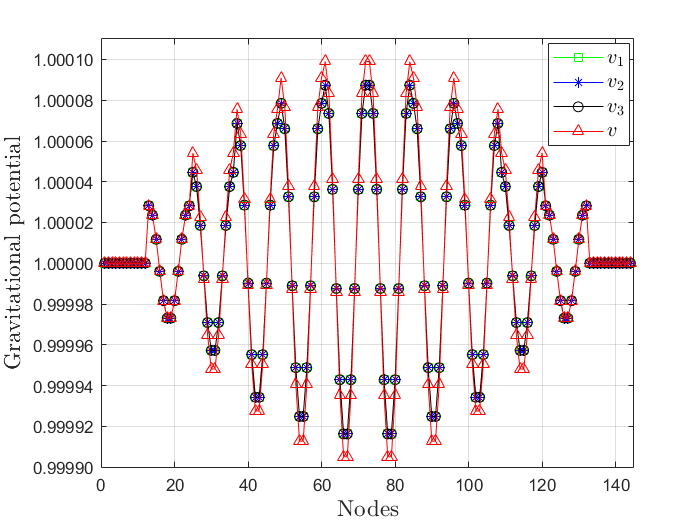}
    \subcaption*{(b)}
     \end{minipage}
    \caption{ (a) Gravitational field approximation $\nabla v_3$ at  144 nodes of mixed (quad/triangular) element meshing of surface of unit sphere and (b) comparison of the first, second and third approximation of gravitational potential with exact solution at 144 mesh nodes for  $\varepsilon =10^{-4}$ for the unit sphere.}
    \label{fig7}
\end{figure}

In this section, the effect of denser meshing of the surface of the unit sphere on the approximate solutions will be analyzed.
 In Fig. \ref{fig6},  a $6\times 6$ meshgrid for $(\phi,\theta)$ has been considered   which provides 36 nodes and 
in Fig. \ref{fig7}, a   $12\times 12$ meshgrid is considered which provides 144 nodes. Such type of meshing is common in spherical and  cylindrical systems which results in a mixed (quad/triangular) element meshing in MATLAB for which 
$v_1$, $v_2$ and $v_3$
are computed. Atkinson \cite{Atkinson} pointed out the disadvantage of this subdivision i.e. it results in a very nonuniform
distribution of nodes and elements where usually there are relatively more nodes
near the poles. Here,
the computed gravitational field vectors $(\nabla v_3)$ for 36 nodes are shown in Fig. \ref{fig6}(a) and for 144 nodes in Fig. \ref{fig7}(a) where they point towards the origin.
A comparison of the approximate solutions of gravitational potential with the exact solution is shown  in Fig. \ref{fig6}(b) and Fig. \ref{fig7}(b) at the corresponding 36 nodes and 144 nodes.   Figs. \ref{fig6}(b) and \ref{fig7}(b) reveal that  near the north and south pole, the gravitational potential approximations almost overlap with the exact solution where the mesh is denser.

 In this case also, the singularity  occurs for different number of  meshing partition, but it has been dealt with the numerical integration  described in \eqref{eq22}. On considering different types of quadrature formulas, singularities  occur  for the same considered  mesh. For example, using MATLAB command \verb+`integral2'+  
which is the double integration with adaptive quadrature formula,  singularities would result  for $6\times  6$ and $12\times 12$ meshgrid. Hence, finding approximation solutions at the surface of the unit sphere strongly depends on both meshgrid partition as well as Gauss quadrature.

However, this problem  seems to persist only in the unit sphere. For $r>1$, any number of meshgrid partition along with quadrature formula provides good results.  This also suggests that the singularity due to collinearity with the origin in the Green's function for the unit sphere  is responsible for the conundrum.
Moreover,
a triangular element meshing of the unit sphere through ABAQUS simulation software
has also been attempted. It has been found that such meshing fails to provide approximation results for all mesh nodes  due to collinearity in the Green's function   during the numerical surface integration on the unit sphere. In  Section 3.2, for $r>1$ with removable singularity, approximation results using triangular element meshing have  been provided in Fig. \ref{fig9}.
Hence, it can be concluded that a mixed (quad/triangular) element meshing  of the surface of the unit sphere shown in Figs. \ref{fig6}(a) and \ref{fig7}(a)  can be used to find the approximation solutions at multiple nodes.

It has been found that an icosahedron as well as mixed element meshing  can provide approximation results for the unit sphere where the quadrature method plays a significant role. 
Hence,
in view of this, an icosahedron is still an efficient way to finding the solution, since computing  approximate solutions at those 12 vertices provides sufficient information on the boundary of unit sphere and hence, earth's  gravitational field.

\subsection{Approximation solutions using perturbation scheme for sphere of radius $r>1$}

\noindent
From the previous sections, it can be noted that
although the approximate solutions are becoming closer to the exact solution, not all the vertices exactly converge to the exact solution for the unit sphere. 
The quadrature method for the numerical integration has to be carefully applied in order to avoid singularities of the Green's function. However,  for $r>1$, the  Green's function possesses removable singularity
(provided in Appendix). Hence, for $r>1$, the perturbation scheme should  provide better approximations.

\begin{figure}[htp]
    \centering
       \begin{minipage}{0.4\textwidth}
 \includegraphics[scale=0.35]{
 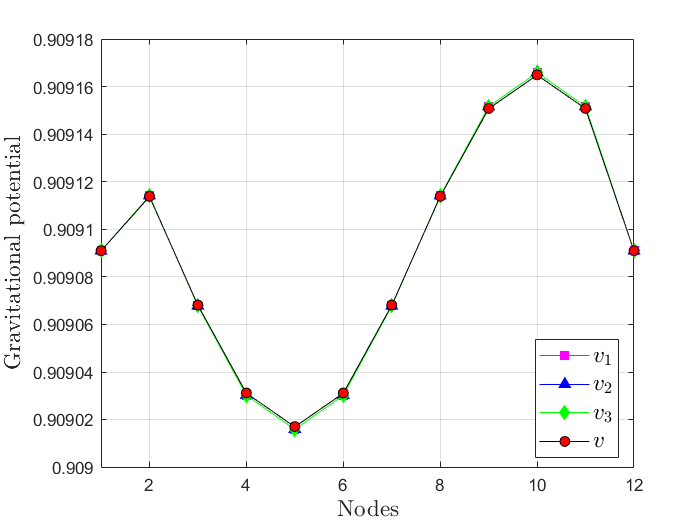}
    \subcaption*{(a)}
    \end{minipage}
    \hspace{1.2cm}
     \begin{minipage}{0.4\textwidth}
    \includegraphics[scale=0.35]{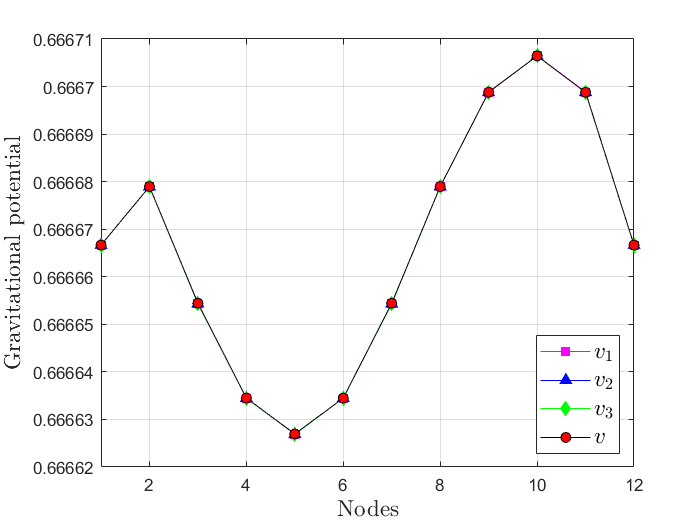}
    \subcaption*{(b)}
     \end{minipage}
    \caption{Graphical representation of the first, second and third approximations and exact solution against vertices of the icosahedron inscribed inside a sphere of radius (a) $r=1.1$ and (b) $r=1.5$ with $\varepsilon =10^{-4}$.}
    \label{fig8}
\end{figure}

\begin{table}[h]
\centering
\begin{tabular}{|c|c|}
\hline
\boldsymbol{ $r$ }  & \boldsymbol{ $|v(P_5) - v_3(P_5) |$  }     \\ \hline
1   & 0.5224208536103170
 $\times 10^{-5}$  \\ \hline
10  & 0.133290600778935  $\times 10^{-11}$ \\ \hline
$10^2$ & 0.133226762955019 $\times 10^{-14}$    \\ \hline
$10^3$ & 0.151788305872824 $\times 10^{-{17}}$    \\ \hline
$10^4$ &  0.135525271663882 $\times 10^{-{19}}$    \\ \hline
$10^5$ & 0.338956580156032  $\times 10^{-{20}}$   \\ \hline
\end{tabular}
\caption{Error analysis of third approximation solution with exact solution at the point $P_5$ for an icosahedron inscribed inside a  sphere of different radius for $\varepsilon=10^{-4}$.}
\label{tab3}
\end{table}

\begin{figure}
    \centering
       \begin{minipage}{0.4\textwidth}
 \includegraphics[scale=0.43]{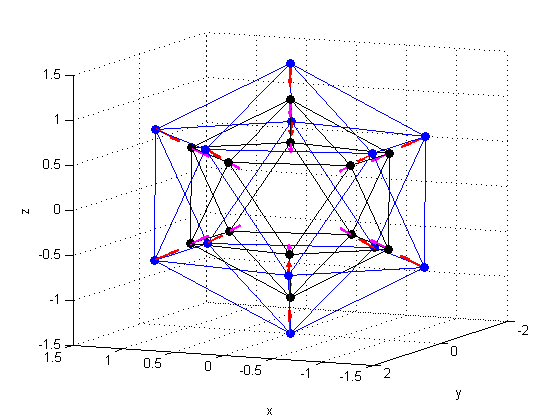}
    \subcaption*{(a)}
    \end{minipage}
    \hspace{1cm}
     \begin{minipage}{0.4\textwidth}
    \includegraphics[scale=0.43]{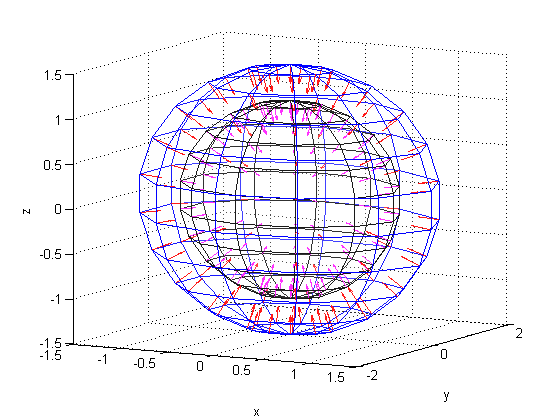}
    \subcaption*{(b)}
     \end{minipage}
    \caption{Gravitational field approximation $\nabla v_3$ at the vertices of (a) two inscribed  icosahedrons with radii $r=1.1$ and $r=1.5$ and (b) mixed element mesh (144 nodes) of two spherical surfaces with radii $r=1.1$ and $r=1.5$.  }
    \label{fig9}
\end{figure}

\begin{figure}
    \centering
       \begin{minipage}{0.45\textwidth}
 \includegraphics[scale=0.4]{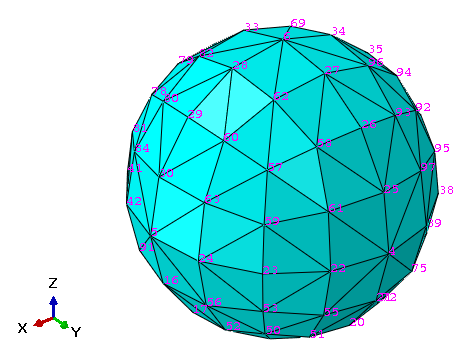}
    \subcaption*{(a)}
    \end{minipage}
    \hspace{-1cm}
     \begin{minipage}{0.4\textwidth}
    \includegraphics[scale=0.38]{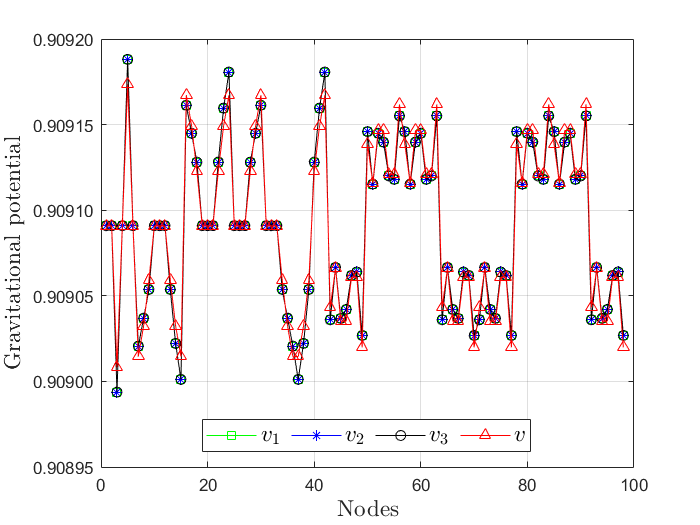}
    \subcaption*{(b)}
     \end{minipage}
      \begin{minipage}{0.4\textwidth}
    \includegraphics[scale=0.38]{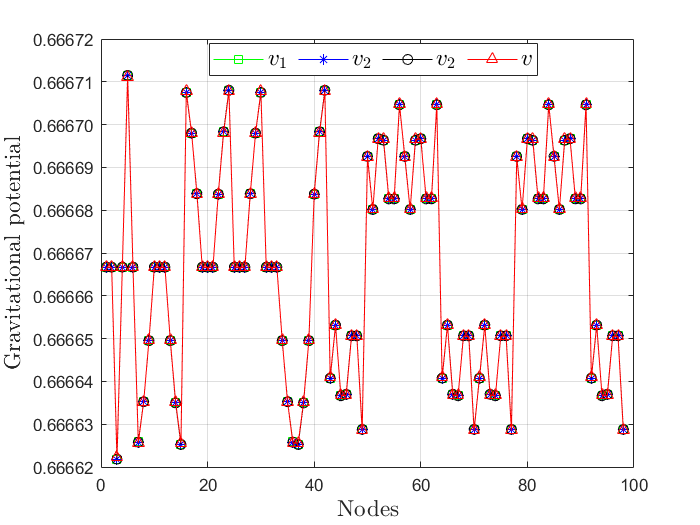}
    \subcaption*{(c)}
     \end{minipage}
     \hspace{1.5cm}
           \begin{minipage}{0.4\textwidth}
    \includegraphics[scale=0.42]{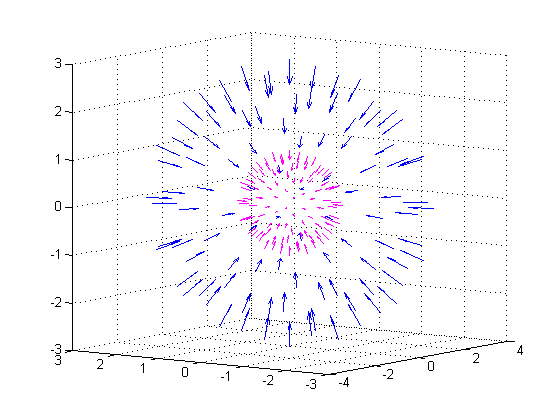}
    \subcaption*{(d)}
     \end{minipage}
    \caption{ (a) Triangular element meshing of the surface of the unit sphere with 98 nodes; (b) comparison of the first, second and third approximation of gravitational potential with exact solution at mesh nodes  for  sphere of  radius $r=1.1$; (c) comparison of the first, second and third approximations of gravitational potential with exact solution at mesh nodes  for   sphere of radius $r=1.5$; and (d) gravitational field approximation $\nabla v_3$ at the vertices of triangular element mesh.
    }
    \label{fig10}
\end{figure}

In order to study the effect of the radius of the sphere, the perturbation scheme for the icosahedron inscribed in a sphere having radius $r=1.1$ and $r=1.5$ is applied. In Fig. \ref{fig8}(a) and \ref{fig8}(b),  
a graphical representation of the approximation solutions ($v_1$, $v_2$ and $v_3$) of gravitational potential along with exact solution ($v$) is depicted for two distinct radii of a sphere. Here, the icosahedron is considered since it has been already proven to be an efficient way of computing the approximation results.  
It can be immediately noticed that as the radius of the sphere ($r>1$) increases, the approximation converges very well to the exact solution. 
It is also to be noted that due to the removable singularity, the numerical integration in this case can be computed without any difficulty, unlike  the unit sphere, for any type of meshing as well as any type of quadrature formula.

In Table \ref{tab3}, an error analysis of the third approximation of gravitational potential at the point $P_5$ is provided for an icosahedron inscribed inside a  sphere of  radius $r=10^k$, for $k=0,1,2,...,5$. As  $r$ increases, it can be clearly noticed that the error reduces drastically.

In Figs. \ref{fig9}(a), third approximation of Gravitational field vectors $\nabla v_3$ at the vertices of   two inscribed icosahedrons with radii r = 1.1 and r = 1.5 are shown whereas in and \ref{fig9}(b), they are shown for a mixed
element mesh (144 nodes) of two spherical surfaces with radii r = 1.1 and r = 1.5.

In Fig. \ref{fig10}(a), a triangular element meshing of the surface of the sphere having 98 mesh nodes is considered. Using ABAQUS simulation software, two different radii, i.e. $r=1.1$ and $r=1.5$ are studied in Figs. \ref{fig10}(b) and \ref{fig10}(c) respectively by keeping the same number of nodes. In this case also, the approximation  converges to the exact solution as $r$ increases. 
It is to be noted that node number 6 represents the north pole and node number 1 represents the south pole in the ABAQUS triangular element meshing.   

In \ref{fig10}(d), the computed gravitational field 
approximations $\nabla v_3$ at the vertices of triangular element mesh of the surface of sphere with radii $r=1.1$ and $r=1.5$ are shown and they point towards the origin.

\subsection{Approximation solutions  for different exact solution}

In this section,  the perturbation scheme is again applied by considering a different exact solution as follows
\begin{align}
    v(r,\phi,\theta)=\frac{1}{r}-\varepsilon\frac{5}{2r^5}\sin{\theta}\;(7\cos^3{\theta}-3\cos{\theta})\; e^{i\phi},
    \label{eq23}
\end{align}

for which the boundary condition is derived as
$(\nabla v\cdot \nabla v)_S=1+\varepsilon h,
$ where
\begin{align*}
   & h(\phi,\theta)= -25e^{i\phi}\sin{\theta} \;(7\cos^3{\theta}-3\cos{\theta})+\varepsilon\frac{25}{4}\left[ 25 e^{2i\phi}\sin^2{\theta} \;(7\cos^3{\theta}-3\cos{\theta})^2\right.\nonumber\\
    &\hspace{1.5cm}\left. + e^{2i\phi}(3-27\cos^2{\theta}+28\cos^4{\theta})^2- e^{2i\phi}(7\cos^3{\theta}-3\cos{\theta})^2 \right].
\end{align*}

\begin{figure}[htp]
    \centering
       \begin{minipage}{0.4\textwidth}
 \includegraphics[scale=0.35]{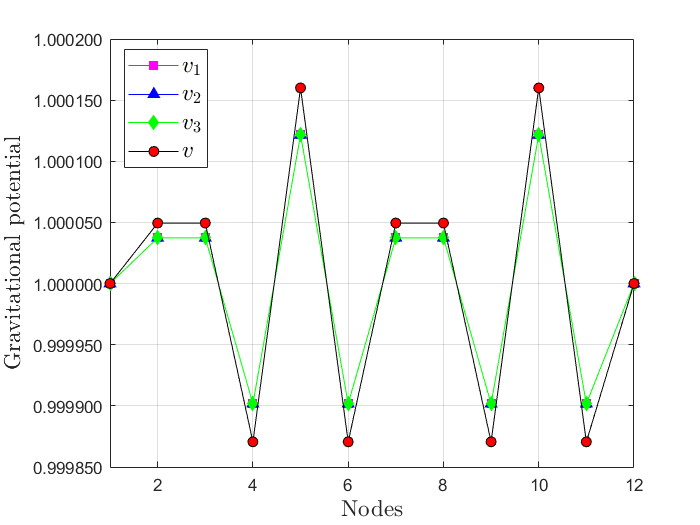}
    \subcaption*{(a)}
    \end{minipage}
    \hspace{1.2cm}
     \begin{minipage}{0.4\textwidth}
    \includegraphics[scale=0.35]{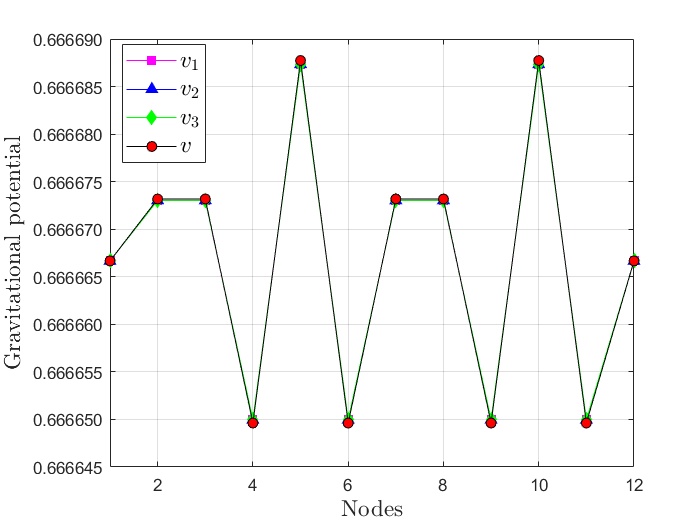}
    \subcaption*{(b)}
     \end{minipage}
    \caption{Graphical representation of the first, second and third approximations along with the new exact solution against vertices of the icosahedron inscribed inside a sphere of radius (a) $r=1$ and (b) $r=1.5$ with $\varepsilon =10^{-4}$.}
    \label{fig11}
\end{figure}

\begin{sidewaystable}[htp]
\scriptsize
\begin{tabular}{|c|c|c|c|}
\hline
$P$   & $(v-v_1) \times 10^{-4}$                            & $(v-v_2) \times 10^{-4}$                            & $(v-v_3) \times 10^{-4}$                            \\ \hline
$P_1$  & 0.000000000000000   + 0.000000000000000i & 0.000000000000000   - 0.000000000000155i & 0.000000000000000   - 0.000000000000155i \\ \hline
$P_2$  & 0.119978006467569 - 0.365998957075105i   & 0.119608639299607 - 0.366266948471041i   & 0.119608612223487 - 0.366266928819236i   \\ \hline
$P_3$  & 0.119978006467569 + 0.365998957074943i   & 0.119608639299607 + 0.366266948472028i   & 0.119608740960508 + 0.366267022622212i   \\ \hline
$P_4$  & -0.312193035505448 - 0.227590007801388i  & -0.312051949867831 - 0.227156388414327i  & -0.312051939513891 - 0.227156420248878i  \\ \hline
$P_5$  & 0.384430058077978 + 0.000000000000000i   & 0.384886621132008 + 0.000000000000909i   & 0.384886495132797 + 0.000000000000909i   \\ \hline
$P_6$  & -0.312193035505448 + 0.227590007801388i  & -0.312051949867831 + 0.227156388413717i  & -0.312051939513891 + 0.227156420248268i  \\ \hline
$P_7$  & 0.119978006467569 - 0.365998957074942i   & 0.119608639301827 - 0.366266948472786i   & 0.119608740962729 - 0.366267022622970i   \\ \hline
$P_8$  & 0.119978006467569 + 0.365998957075106i   & 0.119608639301827 + 0.366266948470903i   & 0.119608612225708 + 0.366266928819098i   \\ \hline
$P_9$  & -0.312193035507669 - 0.227590007803198i  & -0.312051949870051 - 0.227156388416063i  & -0.312051988530238 - 0.227156268962956i  \\ \hline
$P_{10}$ & 0.384430058077978 + 0.000000000000000i   & 0.384886621134228 + 0.000000000000232i   & 0.384886654578587 + 0.000000000000232i   \\ \hline
$P_{11}$ & -0.312193035507669 + 0.227590007803198i  & -0.312051949868941 + 0.227156388416096i  & -0.312051988529127 + 0.227156268962988i  \\ \hline
$P_{12}$ & 0.000000000000000 + 0.000000000000000i   & 0.000000000000000 + 0.000000000000065i   & 0.000000000000000 + 0.000000000000065i   \\ \hline
\end{tabular}
\caption{Error analysis for of the approximation solutions with the new exact solution at the vertices of the icosahedron inscribed inside unit sphere for $\varepsilon=10^{-4}$}.
\label{tab4}
\end{sidewaystable}

In Figs. \ref{fig11}(a) and \ref{fig11}(b), the approximation solutions ($v_1$, $v_2$ and $v_3$) along with the exact solution ($v$) corresponding to  \eqref{eq23} are plotted at the vertices of the icosahedron inscribed in a sphere of radius $r=1$ and $r=1.5$. Similar to the previous case, as  $r$ increases, the approximation solution converges to the exact solution rapidly. However, Fig. \ref{fig11}(a) reveals that  the approximation results contains higher error as compared to what has been found in Fig. \ref{fig4}(a) for the previous considered exact solution  for $r=1$ whereas they both converges rapidly for $r>1$ irrespective of the assumed exact solution. In order to clearly see this, Table \ref{tab4} is provided where it can be noticed that the error for the first, second and third approximation is of order $10^{-4}$ whereas it was $10^{-5}$ for the previous case   (see Table \ref{tab2}). 

Hence, the considered exact solution has a substantial   effect on the approximation results especially at the boundary of the unit sphere ($r=1$). The gravitational fields are not drawn for this case since it will be very similar to Fig. \ref{fig5}.

\newpage

\section{Conclusions}

An important geophysical problem consisted on finding the direction of the external gravitational field of the earth (approximated as a sphere) giving only  its intensity on the surface. The strong nonlinearity of the boundary condition of oblique type placed the problem outside the well known and well studied spectra of boundary value problem on the sphere. Therefore a perturbation scheme was attempted. Suggested by the monopolar nature of the gravitational field, the gravitational potential $v$ of the field was written as the monopole $1/r$ plus an unknown much smaller field $\varepsilon u$, where $\varepsilon$ is a small positive parameter. This perturbation gave rise to a series of exterior Neumann problems on the sphere for each order of the small perturbation parameter.   The Green's function for this problem is known and therefore the solution of the boundary value problem  for each perturbation term was represented analytically as the integral on the sphere of the Green's function times the boundary condition; at each step of the perturbation solution, the boundary condition became more and more complicated since it inherits the strong nonlinear boundary condition.

\noindent
In the present work two major results were achieved:

\begin{enumerate}
    \item[(1)]
 A numerical approximation of the perturbation solution. At this stage the perturbation solution approximates the potential of the gravitational field.

\item[(2)] 
 A numerical approximation of the exterior gravitational field   through a calculation of the gradient of the approximated potential obtained in (1) above.
 
\end{enumerate}

The approximation of the potential of the gravitational field mentioned in (1) was possible using 
the integral representation on the sphere of each term of the solution  in terms of the Green's function. Different meshing (icosahedron, uniform mixed-element meshing  and  dense triangulation) of the surface of the sphere were done  and adaptive quadrature method was applied to calculate the integral. Two examples of exact solutions for the unperturbed  problem were considered in order to be able to check the efficiency of the numerical approximation by the calculation of the error.

In spite of the complexity of the boundary condition as the order of the perturbation increases, the integration of the solution could be carried away without major difficulties. A meshing of the sphere in 12 equally spaced points is induced by taking the vertices of an inscribed icosahedron. This mesh proved to be very efficient to approximate the solutions, i.e.  gravitational potentials and hence gravitational fields, at each perturbation step. Two different types of quadrature formula (Gauss 5-point and adaptive quadrature) depending on the integration parameters, were considered for the surface integration. We call this type of double integration on the icosahedron mesh, the icosahedron method and is an important contribution of this work. This  strategy turns out to be very efficient because it can deal with singularity situation due to collinearity of Green's function. It was also used for denser mixed-element mesh  as well as ABAQUS triangular mesh  in order to compute the approximations of the first three terms of the perturbation scheme. On comparing, it is concluded that the icosahedron method, although simple,  is a much  faster and efficient  way to compute the approximation of the solution with high accuracy. It was also better as compared to applying the numerical surface integral method with applied rotation and denser mesh points. 

Through numerical results   for a sphere of radius $r>1$, it has been concluded that  when the point $P$ is in the exterior domain, a rapid convergence of the approximation solutions to the exact solutions is achieved for all types of meshing, this was expected because the singularities of the Green's function are removable for $r>1$. Furthermore,
in all cases, different values of the small parameter $\varepsilon$ were used and it was observed that as it gets smaller, the approximations get closer to the exact solutions. 
Moreover, approximations of the first three terms were calculated for both examples of exact solutions confirming the above.   
The error in the comparison was of the order of $10^{-5}$ for the first example whereas it was   $10^{-4}$ for the second one which shows the significance of the consideration of exact solution. 

To achieve (2), the numerical approximation of the gravitational field was calculated using the gradient of the computed 
 potentials from (1). Spherical harmonics were used to construct approximation functions at each perturbation step in order to compute the gravitational fields. The approximated field vectors are found to be pointing towards the origin of the unit sphere 
at all mesh nodes (icosahedron, mixed-element mesh, triangular mesh) which matches with gravitational fields of the exact solutions. Further, the  field vectors at the exterior domain of the unit sphere boundary were numerically   computed and  plotted in bigger spherical surfaces outside the unit sphere boundary confirming the radial direction.  
This shows that through the  perturbation scheme, the gravitational fields can be approximately computed at any point on the boundary of the unit sphere as well as in the exterior domain for a given data of surface intensity on the boundary.

\section*{CRediT authorship contribution statement}

{\bf M. S. Chaki:} Conceptualization, Investigation, Methodology, Software, Formal analysis, Writing – original draft.  {\bf Maria  C. Jorge:} Conceptualization,   Formal analysis,  Investigation, Methodology, Writing – original draft, Writing – review \& editing.

\section*{Declaration of competing interest}
The authors declare that they have no known competing financial interests or personal relationships that could have appeared to influence the work reported in this paper.

\section*{Acknowledgement}
MSC would like to thanks Consejo Técnico de la Investigación Científica (CTIC) for providing UNAM Postdoctoral Fellowship. Both the authors acknowledge  PAPIIT DGAPA UNAM (IN101822) project. The authors also thank Dr. Julián Bravo-Castillero for his comments on the manuscript.

\appendix

\setcounter{table}{0}

\renewcommand\thesection{\appendixname}

\end{document}